\documentclass[12pt]{amsart}
\usepackage{amsmath,amssymb}

\newtheorem{thm}{Theorem}[section]

\newtheorem{lemma}[thm]{Lemma}
\newtheorem{cor}[thm]{Corollary}

\newcommand{\R}{{\mathbb R}}

\newcommand{\modo}[1]{{\left|#1\right|}}
\newcommand{\normo}[1]{{\left\|#1\right\|}}
\newcommand{\smodo}[1]{{\mathopen|#1\mathclose|}}
\newcommand{\snormo}[1]{{\mathopen\|#1\mathclose\|}}

\begin{document}
\title[Power-bounded operators]{Power-bounded
operators and related norm estimates}

\author[Kalton]{Nigel Kalton}
\makeatletter
\address{Department of Mathematics,
University of Missouri,
Columbia, MO 65211}
\email{nigel@math.missouri.edu}
\author[Montgomery-Smith]{Stephen Montgomery-Smith}
\makeatletter
\address{Department of Mathematics,
University of Missouri,
Columbia, MO 65211}
\email{stephen@math.missouri.edu}
\urladdr{http://www.math.missouri.edu/\~{}stephen}
\author[Oleszkiewicz]{Krzysztof Oleszkiewicz}
\address{Institute of Mathematics,
Warsaw University,
Banacha 2,
02-097 Warsaw,
Poland}
\email{koles@mimuw.edu.pl}
\author[Tomilov]{Yuri Tomilov}
\address{Department of Mathematics and Informatics,
Nicholas Copernicus University,
Chopin Str. 12/18,
87-100 Torun,
Poland}
\email{tomilov@mat.uni.torun.pl}
\thanks{The first and second named authors were
partially supported
by NSF grants.  The third named author was visiting the
University of Missouri-Columbia
while conducting this research,
and was partially supported by
the Polish KBN Grant 2 P03A 027 22.
The fourth named author was partially supported by the Polish KBN Grant
5 P03A 027 21, and the NASA-NSF Twinning Program.}

\begin{abstract}
\noindent
We consider whether
$ L = \limsup_{n\to\infty} n \snormo{T^{n+1}-T^n} < \infty$ implies
that the operator $T$ is power bounded.
We show that this is so if $L<1/e$, but it
does not necessarily hold if $L=1/e$.
As part of our methods,
we improve a result of Esterle, showing that
if $\sigma(T) = \{1\}$ and $T \ne I$, then
$\liminf_{n\to\infty} n \snormo{T^{n+1}-T^n} \ge 1/e$.
The constant $1/e$ is sharp.
Finally we describe a way to create many generalizations of Esterle's result,
and also give many conditions on an operator which imply that its norm is
equal to its spectral radius.
\end{abstract}

\maketitle

\section{Introduction}

Let $T$ be a bounded linear operator on a complex Banach space $X$.
One of the classical problems in operator theory is to determine
the relation between the size of the resolvent $(T-\lambda I)^{-1}$
when $\lambda$ is
near the spectrum $\sigma(T)$, and the asymptotic properties of
orbits $\{T^n x:  n \ge 0\}$ for each $x \in X$.
The inequality
\[
\snormo{(T-\lambda I)^{-1}} \le \frac{C}{\text{dist}(\lambda, \sigma (T))} ,
\quad \lambda \in \mathbb C \, \backslash \, \sigma (T),
\]
has been extensively studied by,
for example, Benamara and Nikolski \cite{benamara-nikolski}
and also, very recently, by El-Fallah and Ransford \cite{el-fallah-ransford};
see also \cite{lyubich2}, \cite{nagy-zemanek}, \cite{nevanlinna2},
\cite{tomilov-zemanek}.
Such an inequality is extreme in the sense that the converse
inequality (with $C=1$) is always satisfied.
In most cases the relationship to such an inequality and the properties
of the orbits are very difficult to determine.

Thus it is interesting
that one has a very clean equivalence for the resolvent
condition introduced by Ritt \cite{ritt}, which says
there is a constant $C>0$ such that
$$ \snormo{(T-\lambda I)^{-1}} \le \frac C{\modo{\lambda-1}}
   \qquad (\modo\lambda > 1) .$$
Nagy and Zem\'anek \cite{nagy-zemanek}, and independently
Lyubich \cite{lyubich1}, proved the following result
(see also \cite[Theorem 4.5.4]{navanlinna1}).

\begin{thm}
\label{ritt} Let $T$ be an operator on a complex Banach space.
Then $T$ satisfies the Ritt resolvent condition if and only if
\begin{enumerate}
\item $T$ is power bounded, and
\item \label{c T^n+1-T^n}
$\sup_n n\snormo{T^{n+1}-T^n} < \infty$.
\end{enumerate}
\end{thm}

We recall a result of Esterle \cite{esterle} saying that if
$\sigma(T)=\{1\}$ and $T$ is not the identity operator, then
$\liminf_{n\to\infty} n\snormo{T^{n+1}-T^n} \ge 1/12$.  (The
citation given only has $1/96$; this was improved by Berkani
\cite{berkani} to $1/12$.)  Moreover it was noted in \cite[Theorem
4.5.1]{navanlinna1} that if $1$ is a limit point of $\sigma(T)$,
then $\limsup_{n\to\infty} n \snormo{T^{n+1} - T^n} \ge 1/e$. Thus both the
Ritt resolvent condition and condition~(\ref{c T^n+1-T^n}) are
extremal, and it is natural to ask whether these two conditions
are equivalent, at least in the case when $\sigma(T)=1$. Note it
was only recently that Lyubich \cite{lyubich2} constructed
operators satisfying the Ritt condition and $\sigma(T) = \{1\}$.

Another reason that such a question is
interesting is because of the famous Esterle-Katznelson-Tzafriri Theorem
\cite{esterle}, \cite{Ka}, which states that if $T$ is power bounded,
and its spectrum meets the unit circle only at the point $1$, then
$\snormo{T^{n+1}-T^n} \to 0$ as $n \to \infty$.  Thus a positive answer
to our question would provide a partial converse.

Towards this conjecture,
it is known that if $\limsup_{n\to\infty}n\snormo{T^{n+1}-T^n} < 1/12$,
then $T$ is power bounded in a rather trivial manner, that is, it is the
direct sum of an identity operator
and an operator whose spectral radius is less than $1$.
This follows
directly from the result of Esterle cited above.

In this paper, we improve these results.
We answer a conjecture of Esterle \cite{esterle} (see also \cite{berkani})
and show that in his result that $1/12$ may be replaced by $1/e$.
Furthermore an example shows that $1/e$ is sharp.  As a corollary
we show that if $\limsup_{n\to\infty}n\snormo{T^{n+1}-T^n} < 1/e$, then
$T$ is power bounded.  Again we provide an example to show that $1/e$ is
sharp.
In particular, the
condition $\sup_n n\snormo{T^{n+1}-T^n} < \infty$ does not necessarily imply
that $T$ is power bounded.
We leave open the question as to whether it
implies power boundedness in the case that $\sigma(T) = \{1\}$.

Finally we create a general framework which shows how to
easily create results in the same vein as Esterle's result.
For example, one can give conditions concerning $\snormo{T^n-T^m}$
that imply that an operator with $\sigma(T)=\{1\}$ is the identity.
We also give results similar to the special case of Sinclair's
Theorem \cite{sinclair} considered by Bonsall and Crabb \cite{bonsall-crabb},
giving many different
conditions on an operator that imply that its norm is equal
to its spectral radius.

Let us finish this introduction by noting that Blunck \cite{blunck1},
\cite{blunck2} gives
many applications of
the condition $\sup_n n \snormo{T^{n+1}-T^n} < \infty$
to maximal regularity problems.
Also, after this present article was finished, the authors learned of
recent papers \cite{berkani-esterle-mokhtari} and \cite{esterle-mokhtari}
which use similar methods.

Throughout this paper, we will
take the Fourier transform to be
$\hat f(\xi) = \int_{-\infty}^\infty f(x) e^{-i x \xi} \, dx$
and the inverse Fourier transform to be
$\check g(x) = \frac1{2\pi} \int_{-\infty}^\infty g(\xi) e^{i x \xi} \, d\xi$.
All Banach spaces will be complex in the remainder of the paper.

\section{Esterle's Result}

To illustrate the ideas, let us first give a
continuous time version.
The methods used are similar to those in a paper by
Bonsall and Crabb \cite{bonsall-crabb} in their proof of a special case of
Sinclair's Theorem \cite{sinclair}.
The function $W$ described below is often called
the Lambert function (see \cite{corless et al}).

\begin{thm}
\label{esterle continuous}
Let $A$ be a bounded operator on a Banach space such that
$\sigma(A) = \{0\}$.  For each $t>0$ such that
$\snormo{A e^{tA}} \le 1/et$, we have that $\snormo{A} \le 1/t$.
In particular, if
$\liminf_{t\to\infty} t \snormo{Ae^{tA}} < 1/e$,
then $A = 0$.
\end{thm}

\begin{proof}
Let $f(z) = z e^z$.
There is analytic function $W$ such that
$W(f(z)) = z$ in some neighborhood of $0$.  In particular,
by the Riesz-Dunford functional calculus,
$W(tA e^{tA}) = tA$.
Now
$$ W(z) = \sum_{m=1}^\infty p_{m} z^m $$
where, by Lagrange's inversion formula \cite[Ch.~5, Ex.~33]{asmar},
$$ p_{m} = \frac1{m!}\frac{d^{m-1}}{dz^{m-1}} \left(\frac z{f(z)}\right)^m
            \Bigg |_{z=0}
          = \frac{(-m)^{m-1}}{m!} .$$
The radius of convergence of $W$ is $1/e$, and
$ \sum_{m=1}^\infty \modo{p_{m}} e^{-m} = 1 $,
since $f(-1) = -1/e$.
Therefore
$\snormo{W(tA e^{tA})} \le 1$, and the result follows.
\end{proof}

\begin{thm}
\label{esterle}
Let $T$ be a bounded operator on a Banach space such that
$\sigma(T) = \{1\}$.  For each positive integer $n$ such that
$\snormo{T^{n+1}-T^n} \le n^n/(n+1)^{n+1}$, we have that
$\snormo{T-I} \le 1/(n+1)$.  In particular, if
$\liminf_{n\to\infty} n \snormo{T^{n+1}-T^n} < 1/e$,
then $T = I$.
\end{thm}

\begin{proof}
Let $f_n(z) = z(1+z/n)^n$.
There is analytic function $W_n$ such that
$W_n(f_n(z)) = z$ in some neighborhood of $0$.  In particular,
by the Riesz-Dunford functional calculus,
$W_n(n(T^{n+1}-T^n)) = n(T-I)$.
Now
$$ W_n(z) = \sum_{m=1}^\infty p_{nm} z^m $$
where
$$ p_{nm} = \frac1{m!}\frac{d^{m-1}}{dz^{m-1}} \left(\frac z{f_n(z)}\right)^m
            \Bigg |_{z=0}
          = \frac{(-1)^{m-1}}{n^{m-1} (nm+m-1)}
            \binom{nm+m-1}m .$$
The radius of convergence of $W_n$ is $r_n = (n/(n+1))^{n+1}$, and
$\sum_{m=1}^\infty \modo{p_{nm}} r_n^m = n/(n+1)$,
since $f_n(-n/(n+1)) = -r_n$.
Therefore
$\snormo{W_n(n(T^{n+1}-T^n))} \le n/(n+1)$
and the result follows.
\end{proof}

In Section~\ref{general} below, we will generalize this approach
and give many extensions of these results.

Now let us turn out attention to whether the constant $1/e$ in
Theorems~\ref{esterle continuous} and~\ref{esterle} can be
improved. By the results of Lyubich \cite{lyubich2} combined with
Theorem~\ref{ritt}, we know that there must be some upper bound on
the numbers $C>0$ such that $\sigma(T) = \{1\}$ and
$\liminf_{n\to\infty} n \snormo{T^{n+1}-T^n} < C$ imply that
$T=I$. In fact we will be able to modify the examples of Luybich
to show that $C = 1/e$ is sharp.

We will consider the fractional Volterra operators, parameterized
by $\alpha>0$, on $L_p([0,1])$ for $1 \le p \le \infty$, given by the
formula
$$ J^\alpha f(x) = \frac1{\Gamma(\alpha)}\int_0^x (x-y)^{\alpha-1}
   f(y) \, dy ,$$
and also modified fractional Volterra operators
$$ L^\alpha f(x) = \frac1{\Gamma(\alpha)}\int_0^x (x-y)^{\alpha-1}
   e^{y-x} f(y) \, dy .$$
It is well known (and easy to show)
that $(J^\alpha)_{\alpha>0}$ is a $C_0$-semigroup
Similarly $(L^\alpha)_{\alpha>0}$ is also a $C_0$-semigroup.
Thus it is easily seen that $\snormo{(L^\alpha)^n} = \snormo{L^{\alpha n}}
\le 1/\Gamma(\alpha n +1)$, and hence the spectral radius of $L^\alpha$
is zero.

Let us also consider an extension of this operator $\tilde L^\alpha$
on $L_p(\R)$ given by the formula
$$ \tilde L^\alpha f(x)
   = \frac1{\Gamma(\alpha)}\int_{-\infty}^x (x-y)^{\alpha-1}
   e^{y-x} f(y) \, dy .$$
This is a convolution operator.  Therefore,
$\widehat{\tilde L^\alpha f}(\xi) = m_\alpha(\xi) \hat f(\xi)$,
where $m_\alpha$ is the Fourier Transform of $x_+^{\alpha-1} e^{-x}
/\Gamma(\alpha)$.  Direct calculation shows that
$m_\alpha(\xi) = (1+i\xi)^{-\alpha}$,
where here we are taking the principle
branch.

Next,
let $M$ denote the operator of
multiplication by the indicator function of $[0,1]$, then
it is not so hard to see that for any entire function
$f$ we have that $f(L^\alpha) = M f(\tilde L^\alpha) M$, and so
$\snormo{f(L^\alpha)} \le \snormo{f(\tilde L^\alpha)}$.

Now we see that
$\widehat{\tilde L^\alpha e^{-t\tilde L^\alpha}f}(\xi) = k(\xi) \hat f(\xi)$,
where $k(\xi) = m_\alpha(\xi) e^{-t m_\alpha(\xi)}$.  If $0<\alpha<1$, then
$\text{Re}(m_\alpha(\xi)) > 0$, and
$\lim_{\xi\to\pm\infty} \text{arg}(m_\alpha(\xi)) = \alpha\pi/2$.  Hence
it is easy to see that
$$ \limsup_{t\to\infty} t\snormo{L^\alpha e^{-t L \alpha}} \le
\limsup_{t\to\infty} t\snormo{\tilde L^\alpha e^{-t\tilde L^\alpha}} \le
1/e \cos(\alpha \pi/2) .$$

This is enough to show that the constant $C = 1/e$ is sharp in
Theorem~\ref{esterle continuous}.  However, we can do a little better.

\begin{thm}
\label{esterle eg}
\begin{enumerate}
\item
There exists an operator $A \ne 0$ on a Hilbert space, with $\sigma(A)=\{0\}$,
and $\limsup_{t\to\infty} t \snormo{A e^{tA}} \le 1/e$.
\item
There exists an operator $T \ne I$ on a Hilbert space, with $\sigma(T)=\{1\}$,
and $\limsup_{n\to\infty} n \snormo{T^{n+1}-T^n} \le 1/e$.
\end{enumerate}
\end{thm}

\begin{proof}
Let us consider the operator on $L_2([0,1])$
$$ A = -\int_0^{1/2} L^\alpha \, d\alpha .$$
Lyubich \cite{lyubich2} showed that the operator
$B = \int_0^\infty J^\alpha \, d\alpha$ has spectral radius equal to $0$
on $L_p([0,1])$ for all $1 \le p \le \infty$.
Now both $-A$ and $B$ are operators with positive kernels, and the kernel
of $-A$ is bounded above by the kernel of $B$.  It follows that on
$L_p([0,1])$ for $p = 1$ or $p = \infty$ that
$\snormo{A^n} \le \snormo{B^n}$ for all positive integers $n$.
Thus $A$ has spectral radius equal to $0$ on $L_p([0,1])$ for $p=1$ and
$p=\infty$, and hence, by interpolation, for all $1 \le p \le \infty$.

We also define the operator on $L_2(\R)$
$$ \tilde A = -\int_0^{1/2} \tilde L^\alpha \, d\alpha .$$
Following the above argument,
we see that $\snormo{A e^{tA}} \le \snormo{\tilde A e^{t \tilde A}}$, and
that $\widehat{\tilde A e^{t\tilde A}f}(\xi) = k(\xi) \hat f(\xi)$, where
$$ \modo{k(\xi)} =
   \modo{h(\xi)}
   \exp(- t \text{Re}(h(\xi))) ,$$
and
$$ h(\xi) = \int_0^{1/2} m_\alpha(\xi) \, d\alpha .$$
One sees that $\text{arg}(h(\xi)) \to 0$ as $\xi\to\infty$,
and hence it is an easy matter to see that
$\limsup_{t\to\infty} t\snormo{A e^{tA}} \le 1/e$.

The second example is given by $T = e^A$.
Note that $T \ne I$, because otherwise $A = \log(T) = 0$.
The estimate is easily obtained
since $T^{n+1}-T^n = \int_n^{n+1} A e^{tA} \, dt$.
\end{proof}

\section{Power Boundedness}

\begin{thm}
\label{T^n}
Let $T$ be a bounded operator on a Banach space $X$ such that
$\limsup_{n\to\infty} n \snormo{T^{n+1}-T^n} < 1/e $.
Then $X$ decomposes as the direct sum of two closed $T$-invariant
subspaces such that $T$ is the identity on one of these subspaces, and
the spectral radius of $T$ on the other subspace is strictly less
than $1$.  In particular, $T^n$ converges to a projection.
\end{thm}

\begin{proof}
First note that $\sigma(T)$
must be contained in $\{1\} \cup \{z:\modo z < \alpha\}$
for some $\alpha<1$, otherwise it is easy to see that
limit superior of the spectral radius
of $T^{n+1}-T^n$ is at least $1/e$  (see, for example
\cite[Theorem 4.5.1]{navanlinna1}).
Thus there is a projection $P$ that commutes with $T$ such that
$\sigma(T|_{\text{image}(P)}) = \{1\}$, and
the spectral radius of $T|_{\text{ker}(P)}$ is strictly less than $1$.
The result now follows by applying Theorem~\ref{esterle} to
$T|_{\text{image}(P)}$.
\end{proof}

A very similar proof works also for the following continuous time version.
However, we were also able to produce a different proof of this same result.

\begin{thm}
\label{exp tA}
Let $A$ be a bounded operator on a Banach
space $X$ such that
$ L = \limsup_{t\to\infty} t \snormo{A e^{tA}} < 1/e $.
Then $X$ decomposes as the direct sum of two closed $A$-invariant
subspaces such that $A$ is the zero operator on one of these subspaces, and
on the other subspace the supremum of the real part of the spectrum is
strictly negative.
In particular, $e^{tA}$ converges to a projection.
\end{thm}

\begin{proof}
To illustrate the ideas, let us first prove that $e^{tA}$ converges
in the case that
$L < 1/4$, that is,
there are constants $c<1/4$ and $t_0>0$ such
that $\snormo{A e^{tA}} \le c/t$ for $t\ge t_0$.
It follows that $\snormo{A^2 e^{2tA}} \le c^2/t^2$ for $t \ge t_0$, or
$\snormo{A^2 e^{tA}} \le 4c^2/t^2$ for $t \ge 2t_0$.  Then for $t \ge 2t_0$
we have
$$ \snormo{A e^{tA}}
   =
   \normo{\lim_{\tau\to\infty} \int_t^\tau A^2 e^{sA} \, ds}
   \le \frac{4c^2}t, $$
since ${A e^{\tau A}} \to 0$ as $\tau \to \infty$.
Iterating this process,
we get that
$\snormo{A e^{tA}} \le (4c)^{2^k}/4t$ for $t\ge 2^k t_0$.  To put this
another way, $\snormo{A e^{tA}} \le (4c)^{t/2t_0}/4t$ for $t \ge t_0$.
It follows that
$$ e^{t_1A} - e^{t_2A} = \int_{t_2}^{t_1} A e^{sA} \, ds $$
converges to zero as $t_1,t_2 \to \infty$,
that is, $e^{tA}$ is a Cauchy sequence.  Hence it
converges.

The case when $L<1/e$ is only marginally more complicated.
Again,
there are constants $c<1/e$ and $t_0>0$ such that
$\snormo{A e^{tA}} \le c/t$ for $t\ge t_0$.
For any integer $M \ge 2$ we have
that $\snormo{A^M e^{tA}} \le (cM)^M/t^M$ for $t \ge Mt_0$.
Integrating $(M-1)$ times we obtain that
$$ \snormo{A e^{tA}} \le \frac{(cM)^M}{t(M-1)!}  \quad \text{for $t\ge Mt_0$}.$$
A simple computation shows that
$$ \frac{(cM)^M}{(M-1)!} \le \frac Me (ce)^M ,$$
and hence iterating we obtain that if $t>M^k t_0$ then
$$ \snormo{A e^{tA}} \le
   \left(\frac Me\right)^{-1/(M-1)}
   \left( ce \left(\frac Me\right)^{1/(M-1)} \right)^{M^k} \frac 1t.$$
By choosing $M$ is sufficiently large, we see that
there exist constants $c_1,c_2>1$ such that
$\snormo{A e^{tA}} \le c_1 c_2^{-t}/t$ for $t \ge t_0$,
and hence
$\snormo{e^{tA}}$ converges.

Now it is clear that $S=\lim_{t \to \infty} e^{tA}$ is a
bounded projection (because $S^{2}=S$) such that $Se^{tA}=e^{tA}S=S$.
Let $X_1 = \text{Im}(S)$, and $X_2 = \text{Ker}(S)$, so
$X = X_1 \oplus X_2$.  These spaces are
clearly invariant under $e^{tA}$, and hence
invariant under $A = \lim_{t\to0}(e^{tA}-I)/t$.
Since
$S|_{X_1} = I|_{X_1}$ we see immediately that $e^{tA}|_{X_1} = I|_{X_1}$,
and so $A|_{X_1} = \lim_{t\to0}(e^{tA}|_{X_1}-I|_{X_1})/t = 0$.
Furthermore, we have that $e^{tA}|_{X_2} \to 0$.  Let $t_0$ be such that
$\snormo{e^{t_0 A}|_{X_2}} \le 1/2$.  Then the spectral radius of
$e^{t_0 A}|_{X_2}$ is bounded by $1/2$, and so
$\sup\text{\rm Re}(A|_{X_2}) < -\log(2)/t_0$.
\end{proof}

We also point out that that one could prove Theorem~\ref{T^n} in a similar
manner.  But the details can be quite complicated.  It is also possible
to deduce Theorem~\ref{T^n} from Theorem~\ref{exp tA}.  Briefly,
if $\snormo{T^{n+1}-T^n} \le (1+\epsilon)L/(n+1)$ for large enough enough $n$,
then by writing out the power series for $(T-I)e^{tT}$ about $t=0$
one obtains that
$\snormo{(T-I)e^{tT}} \le (1+2\epsilon)Le^t/t$ for large enough $t$.
The result now follows quickly by applying Theorem~\ref{exp tA} to $A = T-I$,
remembering that $\sigma(T) \subset \{1\} \cup \{z:\modo z < 1\}$.

Now we give some counterexamples to show that in general the condition
$\sup_n n\snormo{T^{n+1}-T^n} < \infty $ does not necessarily imply
power boundedness.

\begin{thm}
There exists a bounded operator $T$ on $L_1(\R)$ such that
$ \sup_n n\snormo{T^{n+1}-T^n} < \infty $,
and $\snormo{T^n} \approx \log n$.
\end{thm}

\begin{proof}
The example is a multiplier on $L_1(\R)$ given by
$\widehat{Tf}(\xi) = m(\xi) \hat f(\xi)$.
It is well known that such an operator is bounded if the inverse
Fourier transform
$\check m$ is a measure of bounded variation, and indeed that the norm
is equal to the variation of $\check m$.

Let us consider the case
$$ m(\xi) = \left\{
   \begin{array}{cl}
     1 & \text{if $\modo\xi \le 1$}\\
     \exp(1-\modo{\xi}) & \text{if $\modo\xi > 1$}.
   \end{array} \right.
$$
An explicit computation
shows that the inverse Fourier transform of $m^n$ is
$$ \frac{n x\cos(x)+n^2\sin(x)}{\pi x(x^2+n^2)} $$
and that the inverse Fourier transform of $m^{n+1}-m^n$ is
$$ \frac
   {(x^2-n(n+1))\cos(x) + (2n x+x)\sin(x)}
   {\pi (x^2+n^2)(x^2+(n+1)^2)} ,$$
and it is now easy to verify the claims.
\end{proof}

\begin{thm}
On any Banach space $X$, there exists a bounded operator $T:X\to X$
such that
$ \limsup_{n\to\infty} n\snormo{T^{n+1}-T^n} < \infty $,
and $\snormo{T^n} \to \infty$.  Furthermore there is an equivalent
norm $\modo{\,\cdot\,}$ on $X$ so that
$ \limsup_{n\to\infty} n\smodo{T^{n+1}-T^n} \le 1/e $.
\end{thm}

\begin{proof}
In any Banach space $X$ we may find a sequence $e_n\in X$ with
$\|e_n\|=1$ and bi-orthogonal functionals $e_n^*\in X^*$ such that
$\sup_n\|e_n^*\|=M<\infty$ and such that $(e_n)_{n=1}^{\infty}$ is
not a basic sequence.  Indeed, by \cite{PS}, any subspace of $X$
with a basis has a normalized conditional basis, which may be
re-ordered to give the example.  (We remark that if $X$ is
separable, then one can choose $(e_n)_{n=1}^{\infty}$ to be
fundamental by using \cite{OP} or \cite{P}).  We refer to
\cite{LT} for details.

Let $E=[e_n]_{n=1}^{\infty}$ be the closed linear span of
$(e_n)_{n=1}^{\infty}$.  Define $T:X\to X$ by
$$Tx = x+\sum_{k=1}^{\infty}(\lambda_k-1)e_k^*(x)e_k$$
where $\lambda_k=\exp(-1/k!)$.  Since $|\lambda_k-1|\le 1/k!$ it
follows that $T$ is bounded and $\|T\|\le e+1$.

Consider
$$(T^n-T^{n+1})x=
\sum_{k=1}^{\infty}(\lambda_k^n-\lambda_k^{n+1})e_k^*(x)e_k.$$
Hence
$$ n\|T^n-T^{n+1}\|\le M\sum_{k=1}^{\infty}\frac{ ne^{-n/k!}}{k!}.$$
To estimate this sum suppose $m!<n\le (m+1)!$.  Then
$$\sum_{k=1}^{\infty}\frac{ ne^{-n/k!}}{k!}=
\left(\sum_{k=1}^{m-1} \frac{n}{k!}e^{-n/k!}\right) +
\frac{n}{m!}e^{-n/m!} +
\left(\sum_{k=m+1}^{\infty}\frac{n}{k!}e^{-n/k!}\right).$$ Simple estimates
show that the two sums converge to $0$ as $n \to \infty$, and it
is easy to see that the middle term is bounded by $1/e$. Hence
$\limsup_n n\|T^n-T^{n+1}\| \le M/e$.

Now we claim that if $\sup\|T^n\|<\infty$ then $(e_n)$ is a basic
sequence, giving a contradiction. To do this we estimate $\|P_n\|$
where
$$ P_nx=\sum_{k=1}^n e_k^*(x)e_k.$$
Then
$$ P_nx+T^{n!}x= x + \sum_{k=1}^n \lambda_k^{n!}e_k^*(x)e_k
+\sum_{k=n+1}^{\infty}(\lambda_k^{n!}-1)e_k^*(x)e_k.$$ Thus
$$ \|P_n+T^{n!}-I\| \le \sum_{k=1}^n e^{-n!/k!} +
\sum_{k=n+1}^{\infty} \frac{n!}{k!}.$$ As before we can estimate
both sums to be uniformly bounded in $n$. So if $T$ is power-bounded
then $(P_n)$ is uniformly bounded, and hence
$(e_n)_{n=1}^{\infty}$ is basic.

Let us remark that the above construction also yields a
counter-example if $X$ is reflexive and $(e_n)_{n=1}^{\infty}$ is
a basis of an uncomplemented subspace of $X$, since in that case
one can show that $P_n$ converges in the weak-operator topology to
a projection on $E$.

To obtain the equivalent norm on
$X$, set $|x|=\max(\|x\|,\sup_n|e_n^*(x)|)$.  Let $X=(X,|\cdot|)$
and note that in this case $M=1$.
\end{proof}

\section{A general approach}
\label{general}

In this section we will discuss how to extend Theorems
\ref{esterle continuous} and \ref{esterle} by a more general
approach. We first isolate the argument used.

To do this, let us introduce a class of analytic functions.  Let
$f$ be an analytic function defined on a disk $\{z:|z|<R\}$ (we
allow the case when $f$ is entire and $R=\infty$).

We will say
that $f\in \mathcal P$ if:\begin{enumerate}
\item $f(0)=0$.\item $f'(0)\neq 0$.\item $f(x)\in\mathbb R$ if
$-R<x<R$.
\item The local inverse function $\varphi=f^{-1}$ of $f$ at the origin,
which is defined in a
neighborhood of $0$ with $\varphi(0)=0$, satisfies the conditions
$\varphi^{(n)}(0)\ge 0$ for all $n\ge 1$.
\end{enumerate}

We remark that in \cite{bonsall-crabb} the key idea is that
$f(z)=\sin z$ is in class $\mathcal P$.  In \S2, we essentially
used the fact that the functions $ze^{-z}$ and $z(1-\frac{z}n)^n$
are in class $\mathcal P$. Before proceeding let us include
another simple example which illustrates the basic ideas. During
the late 1960's a series of papers investigated conditions on the
sequence of norms $\|I-T^n\|$ which imply that $T=I$.  A typical
result is that of Chernoff \cite{chernoff}, that says if
$\sup_{n\ge 0}\|I-T^{2^n}\|<1$ then $T=I$. 
Later Gorin \cite{gorin} considered
similar results for sequences $(q_n)_{n=0}^{\infty}$ replacing
$(2^n)$; he showed the result is also true for sequences
 $q_n=3^n,4^n,5^n$ but not $6^n$. More generally the conclusion 
is true if $q_0=1$ and $q_{n+1}/q_n\le 5$.
Let us prove the following
simple result:

\begin{thm}\label{chernoffgorin} Suppose $T$ is a bounded operator on a
Banach space
$X$.  Suppose $\lambda=1$ is the only complex
solution of the
system of inequalities
$$ |1-\lambda^n|\le \|I-T^n\| \qquad n=1,2,\ldots$$
Then $T=I$.\end{thm}

\begin{proof} It is clear that $\sigma(T)=\{1\}$.
Assume $0<a<1$.  Then there exists $n\in\mathbb N$ so that
$\|I-T^n\|<1-a^n$. Consider the function $f(z)=1-(1-z)^n$.  This
is in class $\mathcal P$ and $\varphi$ is given by $\varphi(z)=
1-(1-z)^{\frac1n}$ for $|z|<1$. Let $A=I-T$ so that $A,f(A)$ are
quasi-nilpotent.  By the Riesz-Dunford functional calculus $$
A=\varphi(f(A))=\sum_{k=0}^{\infty}\frac{\varphi^{(k)}(0)}{k!}f(A)^k.$$
In particular $\|A\| \le \varphi(\|f(A)\|)<1-a$.  It follows that
$A=0$ and $T=I$.\end{proof}

We now derive a Corollary which is a slightly stronger form of the
results of Gorin cited above.  Note that if $c<5$ we have
$2\sin({\pi}/{(c+1)})>1$.

\begin{cor}\label{chernoffgorin2} Suppose $T$ is an operator on a Banach
space such that $\liminf_{n\to\infty}\|I-T^n\|<1$.  Suppose  for
some $c>1$ there is a sequence $(q_n)_{n=0}^{\infty}$ with $q_0=1$
and $q_{n+1}\le cq_n$  if $n\ge 0$ such that $\|I-T^{q_n}\|<
2\sin({\pi}/{(c+1)})$ for $n\ge 0$.  Then $T=I$.\end{cor}

\begin{proof} Both statements follow very simply from the Theorem.
 Indeed if $|1-\lambda^n|\le \|I-T^n\|$ for all $n$ then the fact
 that $\liminf_{n\to\infty}\|I-T^n\|<1$ is enough to imply
 $|\lambda|=1$.  Now if $\lambda=e^{i\theta}$ where $|\theta|\le \pi$ we have
 $|\theta|<2\pi/(c+1)$.  If $\theta\neq 0$ let $N$ be the least integer such that
 $q_{N+1}|\theta|\ge 2\pi/(c+1)$.  Then $q_{N+1}|\theta|\le cq_N|\theta|\le 2c\pi/(c+1)$ so that
 $|1-\lambda^{q_{N+1}}|\ge 2\sin(\pi/(c+1))$.  This yields a contradiction and so $\lambda=1$.
\end{proof}

Our next Lemma gives us a recipe for constructing next examples of
functions in class $\mathcal P$, when explicit calculation of the
inverse function $\varphi$ may be difficult.

\begin{lemma}\label{newfunctions}  Let $f,h$ be analytic functions on the disk $\{z:\
|z|<R\}$.  Suppose $f\in\mathcal P$ and that
  $h$ satisfies $h(0)>0$, $h^{(n)}(0)\ge 0$ for all $n\ge 1$ and $h$ is
nonvanishing.  Then if $F(z)=f(z)/h(z)$ we have $F\in\mathcal
P$.\end{lemma}

\begin{proof} The first three conditions are obvious.   For the last condition,
let $\varphi$ be the local inverse of $f$ at the origin defined on
some disk centered at the origin.  Let $0<\rho<\frac12$ be chosen
so that $\rho$ is smaller than the radius of convergence of the
power series expansions of $h$ and $\varphi$ around the origin and
let $M\ge 1$ be an upper bound for $|h|,|h'|,|\varphi|$ and
$|\varphi'|$ on the disk $ \{z:\ |z|\le \rho\}$.   For fixed $w$
consider the map $\Phi_w(z)= \varphi(wh(z))$ for $|z|\le \rho$.
Then if $M|w|<\rho$, we have $|\Phi_w(z)|\le M|w||h(z)|\le
M^2|w|$.  Thus if $|w|<M^{-2}\rho$ we have that $\Phi_w$ maps
$\{z:\ z\le \rho\}$ to itself. We also have $ |\Phi_w'(z)|\le
M^2|w|<\rho$. We conclude that if $|w|<M^{-2}\rho$ then $\Phi_w$
maps the disk $\{z:\ |z|\le \rho\}$ to itself and satisfies
$|\Phi_w'(z)|\le \frac12$ for $|z|\le \rho$. By the Banach
contraction mapping principle if $|w|<M^{-2}\rho$ we can define
$g_n(w)$ by $g_n(0)=0$ and then $g_n(w)=\Phi_w(g_{n-1}(w))$ and
$g_n(w)$ converges to the unique fixed point $\psi(w)$ of
$\Phi_w$. The convergence is uniform on the disk $\{w:\
|w|<M^{-2}\rho\}$.  By induction each $g_n$ is analytic and has
non-negative coefficients in its Taylor series expansion about the
origin.  It follows that $\psi$ has the same properties, and
$\psi$ is clearly the inverse function of $F$.\end{proof}

Let us say $f\in\mathcal P$ is {\it admissible} if there exists
$0<x<R$ such that $f'(x)=0$. If $f$ is admissible let $\xi$ be the
least positive solution of $f'(x)=0$ and suppose $\delta$ is the
radius of convergence of the power series expansion of $\varphi$.

\begin{lemma}\label{admissible} If $f$ is admissible then
$\delta=f(\xi)$ and $$\xi=
\sum_{k=0}^{\infty}\frac{\varphi^{(k)}(0)}{k!}f(\xi)^k.$$
\end{lemma}

\begin{proof} Clearly we have $\varphi(x)<\xi$ if $0<x<\delta$.  Let
$\eta=\lim_{x\to \delta}\varphi(x)$ so that $\eta\le \xi$. If
$\eta=\xi$ we are done. Assume $\eta<\xi$. Then it is clear that
$\varphi'$ is bounded above by $L=f'(\eta)^{-1}$. Let
$U=\{\varphi(z):\ |z|<\delta\})$. Let $U_n=\{z:\
d(z,U)<\frac1n\}$. Then $U$ is contained in the disk $\{z:\
|z|<\eta\}$ and so for large enough $n$, $U_n$ is contained in the
domain of $f$.  Then $f$ cannot be univalent on any $U_n$, for, if
it were, $\varphi$ could be extended to an analytic function on a
disk of radius greater than $\delta$.  Pick $z_n,w_n\in U_n$ so
that $w_n\neq z_n$ and $f(w_n)=f(z_n)$.   We can find $w,z\in
\overline{U}$ so that $(w,z)$ is an accumulation point of
$(w_n,z_n)$.  If $w=z$ then $f'(w)=0$ and this implies $\varphi'$
cannot be bounded above, yielding a contradiction.  If $w\neq z$
then we choose $u_n,v_n$ with $|u_n|<r, |v_n|<r$ and
$\varphi(u_n)\to w,\ \varphi(v_n)\to z$.  Then $u_n,v_n\to
f(w)=f(z)$ but
$$ |w-z|\le \limsup_{n\to\infty}L|u_n-v_n|=0.$$  This also yields
a contradiction and the proof is complete.\end{proof}

\begin{thm}\label{nilpotent}  Let $A$ be a quasi-nilpotent operator on a Banach
space $X$.  Suppose $f$ is an admissible analytic function defined
on a disk $\{z:\ |z|<R\}$ and suppose $\xi$ is the smallest
positive solution of $f'(x)=0$.  Then if $\|f(A)\|<f(\xi)$ we have
$\|A\|< \xi$.\end{thm}

\begin{proof} Let $\varphi$ be the local inverse at the origin.
Then we have
$$
A=\varphi(f(A))=\sum_{n=0}^{\infty}\frac{\varphi^{(n)}(0)}{n!}(f(A))^n.$$
Hence by Lemma \ref{admissible}
$$ \|A\|<
\sum_{n=0}^{\infty}\frac{\varphi^{(n)}(0)}{n!}f(\xi)^n=\xi.$$
\end{proof}

Let us note that at this point that we can recapture Theorems
\ref{esterle continuous} and \ref{esterle} (without computing
derivatives explicitly).  Indeed $z$ belongs to $\mathcal P$ and
hence $f(z)=ze^{-z}$ is admissible with $\xi=1$ and $f(\xi)=1/e$.
Similarly $ f(z)=(1-z)^n-(1-z)^{n+1}=z(1-z)^n$ is admissible with
$\xi=1/(n+1)$ and $f(\xi)=n^n(n+1)^{-n-1}$.

Let us now extend these results slightly.  The first theorem below
is a trivial application of the same ideas.

\begin{thm} Suppose $A$ is a quasi-nilpotent operator and for some
positive integer $m$, $\|Ae^{-A^m}\|< (me)^{-1/m}$.  Then $\|A\|<
m^{-1/m}$.  Hence if
$\liminf_{t\to\infty}\|tAe^{-t^mA^m}\|<(me)^{-1/m}$ then $A=0$.
\end{thm}

\begin{thm}\label{esterle2} Suppose $T$ is a bounded operator with
$\sigma(T)=\{1\}$ and for some $m>n\in\mathbb N$ we have
$$\|T^m-T^n\|<\left(1-\frac{n}{m}\right) \left(\frac{n}{m}\right)^{n/(m-n)}.$$  Then
$\|T-I\|<1-(\frac{n}{m})^{1/(m-n)}$.\end{thm}

\begin{proof} We show that $f(z)=(1-z)^n-(1-z)^m$ is admissible.
This follows from Lemma \ref{newfunctions} since $f(z)= (1-z)^n
(1-(1-z)^{m-n})$ and the function $1-(1-z)^{m-n}$ is in $\mathcal
P$ since its local inverse at the origin is given by
$1-(1-z)^{1/(m-n)}$.  Now apply Theorem \ref{nilpotent} to $I-T$.
\end{proof}

It is possible to derive  other formulas of the type of Theorem
\ref{esterle} from Theorem \ref{esterle2}.  For example we have
the following Corollaries:

\begin{cor}  Suppose $T$ is a bounded operator with
$\sigma(T)=\{1\}$. If $$ \liminf_{m/n\to\infty}\|T^m-T^n\|<1$$
then $T=I$.

More precisely if
$$ \limsup_{m/n\to\infty} \frac{m}{n\log(m/n)}(1-\|T^m-T^n\|)>1$$
then $T=I$.\end{cor}

\begin{cor}
Suppose $T$ is a bounded operator with $\sigma(T)=\{1\}$. If $$
\liminf_{p/n\to 0}\frac{n}{p}\|T^{n+p}-T^n\|<\frac1e$$ then $T=I$.
\end{cor}

\begin{cor}
Suppose $T$ is a bounded operator with $\sigma(T)=\{1\}$. Suppose
$0<s<1$. If
$$ \liminf_{\substack{m/n\to s\\ m,n\to\infty}}
\|T^{m}-T^n\|<(1-s)s^{s/(1-s)}$$ then
$T=I$.\end{cor}

The next theorem is a generalization of the argument used
by Bonsall and
Crabb \cite{bonsall-crabb}
to prove a special case of Sinclair's Theorem \cite{sinclair}, namely that
the norm of an hermitian element $A$ of a Banach algebra coincides
with its spectral radius $r(A)$.

\begin{thm}\label{sinclairtype}  Suppose $f$ is an admissible
entire function.  Suppose that for every $-\pi<\theta\le \pi$ we
have either:
\begin{enumerate}
\item
\label{11}
$\displaystyle\sup_{t>0}|f(te^{i\theta})|> f(\xi)$,
or
\item
\label{12}
$\displaystyle|f(te^{i\theta})|<f(\xi)$ for $0<t<\xi $.
\end{enumerate}
Let $A$ be
any operator satisfying
$$ \sup_{t>0}\|f(tA)\|\le f(\xi).$$ Then
$ r(A)=\|A\|$. In particular, if $A$ is quasi-nilpotent then
$A=0$. Furthermore if
$$ \sup_{t>0}\|f(tA)\|< f(\xi)$$ then $A=0$.
\end{thm}

\begin{proof}  We start by observing that if $\lambda\in\sigma(A)$
then $\sup_{t>0} |f(t\lambda)|\le f(\xi)$.  Let $r=r(A)$.  If
$tr<\xi$ then  by (\ref{11}) and (\ref{12}) we have
$|f(t\lambda)|<f(\xi)$ for every $\lambda\in\sigma(A)$.
  Thus applying the Riesz-Dunford functional calculus to $tA$ we have $tA=\varphi (f(tA))$
and so
$$ t\|A\|<\sum_{n=0}^{\infty}\frac{\varphi^{(n)}(0)}{n!}f(\xi)^n=\xi.$$
Hence $\|A\|< \xi/t$ and it follows that $\|A\|\le r(A)$.

For the last part of the theorem, assume that $\sigma(A)\neq
\{0\}$.  Then there exists $-\pi<\theta\le \pi$ with
$\sup_{t>0}|f(te^{i\theta})|<f(\xi)$.  It is easy to see that this
implies that $\varphi$ is unbounded on the disk $\{z:\
|z|<f(\xi)\}$ which contradicts Lemma \ref{admissible}.  Hence $A$
is quasi-nilpotent and the conclusion follows.\end{proof}

In the Bonsall-Crabb argument for Sinclair's theorem one takes
$f(z)=\sin z$ and shows that it verifies the hypotheses and hence
$\|\sin tT\|\le 1$ for all $t>0$ implies that the norm and
spectral radius of $T$ coincide. Other functions are permissible
however, and lead to more general results of this type:

\begin{thm}\label{appl} Let $A$ be an operator on a Banach space $X$.  Then
each of the following conditions implies that $r(A)=\|A\|$.
\begin{enumerate}
\item $\displaystyle\sup_{t>0} t\|Ae^{-tA}\|\le e^{-1}$.
\item $\displaystyle\sup_{t>0} t\|Ae^{-tA^m}\|\le (me)^{-1/m}$ for
$m>1$ an integer.
\item $\displaystyle\sup_{t>0}\|e^{-tA}-e^{-stA}\|\le (s-1) s^{-s/(s-1)}
$ for some $ s>1$.
\item $\displaystyle\sup_{t>0}\|e^{-(s+i)tA}-e^{-(s-i)tA)}\|\le
\frac{2e^{-s\arctan(1/s)}}{\sqrt{1+s^2}}$ for some $s\ge 0$.
\end{enumerate}  In each case a strict inequality implies that
$A=0$.
\end{thm}

\begin{proof} The first two are immediate deductions from the
preceding Theorem \ref{sinclairtype}. We then must show for the
remaining cases that $e^{-z}-e^{-sz}$ for $s>1$ and
$e^{-sz}\sin{z}$ for $s>0$ satisfy the conditions of Theorem
\ref{sinclairtype} (the case $s=0$ is Sinclair's theorem).

Note first that $f(z)=e^{-z}(1-e^{-(s-1)z})$ is admissible by
Lemma \ref{newfunctions}, since $1-e^{(s-1)z}\in\mathcal P$.  In
this case $\xi=(s-1)^{-1}\log s$ and $f(\xi)<1$.  Let us assume
$-\pi<\theta<\pi$ and $\theta\neq 0$.  If $|\theta|>\frac{\pi}2$
then $f(te^{i\theta})$ is unbounded; if $|\theta|=\frac{\pi}2$
then $\sup_{t>0}|f(te^{i\theta})| = 2>1$. If
$|\theta|<\frac{\pi}2$ then we observe that
$$ |f(te^{i\theta})| =
e^{-t\cos\theta}|1-e^{-(s-1)te^{i\theta}}|.$$ Assume that
$\sup_{t>0}|f(te^{i\theta})|\le f(\xi)$. Pick $t_0$ so that
$(s-1)t_0|\sin\theta|=\frac{\pi}2$. Then
$$ e^{-\xi} > f(\xi)\ge |f(t_0e^{i\theta})|\ge e^{-t_0\cos\theta}.$$  Hence
$ t_0\cos\theta>\xi$.  Choose $t_1<t_0$ so that
$t_1\cos\theta=\xi$.  Then $|f(t_1e^{i\theta})|\le f(\xi)$ implies
that $(s-1)t_1|\sin\theta|$ is a multiple of $2\pi$.  Since
$t_1<t_0$ this is impossible.

Next consider $f(z)=e^{-sz}\sin{z}$ where $0<\theta<\frac{\pi}2$.
In this case $\xi=\arctan s^{-1}$.  We can again use Lemma
\ref{newfunctions} to see that $f$ is admissible. Clearly if
$|\theta|\ge\frac{\pi}2$ then $f(te^{i\theta})$ is unbounded on
$\{t>0\}$.  If $0<|\theta|<\frac{\pi}2$ we use the fact that if
$z=x+iy$ then $$ |f(z)|\ge e^{-sx}\cosh y|\sin x|.$$ Hence $
|f(te^{i\theta})|> |f(t\cos\theta)|$ and so
$\sup_{t>0}|f(te^{i\theta})|>f(\xi)$.\end{proof}

\end{document}